\documentclass[12pt]{article}
\usepackage{epic,eepic,epsf,epsfig}
\usepackage{amsmath,enumerate}
\usepackage{amssymb}
\usepackage{amsbsy}

\usepackage{amsfonts}
\newtheorem{theorem}{Theorem}[section]%
\newtheorem{lemma}[theorem]{Lemma}%
\newtheorem{cor}[theorem]{Corollary}%
\newtheorem{remark}[theorem]{Remark}%

 \def\Omega{\Omega}

\def\f{\noindent}

\newcommand{\qed}{\mbox{\raisebox{0.7ex}{\fbox{}}} \vspace{4truemm}}
\def\demo{\f {\bf Proof.}\hskip10pt}

\begin{document}

\baselineskip 16pt

\title{ \vspace{-1.2cm}
Finite groups with some particular maximal invariant subgroups being nilpotent or all non-nilpotent maximal invariant subgroups being normal
\thanks{\scriptsize This research was supported in part by Shandong Provincial Natural Science Foundation, China (ZR2017MA022)
and NSFC (11761079).
\newline
 \hspace*{0.5cm} \scriptsize $^{\ast\ast}$Corresponding
  author.
\newline
       \hspace*{0.5cm} \scriptsize{E-mail addresses:}
       shijt2005@163.com\,(J. Shi),\,\,xufj2023@s.ytu.edu.cn\,(F. Xu).}}

\author{Jiangtao Shi\,$^{\ast\ast}$,\,\,Fanjie Xu\\
\\
{\small School of Mathematics and Information Sciences, Yantai University, Yantai 264005, China}}

\date{ }

\maketitle \vspace{-.8cm}

\begin{abstract}
Let $A$ and $G$ be finite groups such that $A$ acts coprimely on $G$ by automorphisms. We provide a complete classification of a finite group $G$ in which every
maximal $A$-invariant subgroup containing the normalizer of some $A$-invariant Sylow subgroup is nilpotent. Moreover, we show that both the hypothesis that every
maximal $A$-invariant subgroup of $G$ containing the normalizer of some $A$-invariant Sylow subgroup is nilpotent and the hypothesis that every non-nilpotent maximal
$A$-invariant subgroup of $G$ is normal are equivalent.

\medskip \f {\bf Keywords:} maximal invariant subgroup; normalizer; nilpotent; solvable; normal\\
{\bf MSC(2020):} 20D10; 20D25
\end{abstract}

\section{Introduction}\label{s1}

All groups are assumed to be finite. It is known that a non-nilpotent group all of whose
maximal subgroups are nilpotent is called a minimal non-nilpotent group, and R\'{e}dei [7]
provided the structure of a minimal non-nilpotent group. As a generalization of Schmidt's
theorem, Deng, Meng and Lu [4, Theorem 3.1] described the structure of groups of even
order in which all maximal subgroups of even order are nilpotent. For any fixed prime
divisor $p$ of $|G|$, the order of a group $G$, Shi, Li and Shen [9, Theorem 1.8] proved that
if every maximal subgroup of $G$ of order divisible by $p$ is nilpotent then $G$ is either $q$-
nilpotent or $q$-closed for each prime divisor $q$ of $|G|$. Furthermore, Shi and Tian [11, Theorem 1.1] gave a complete classification of a group $G$ in which all maximal subgroups
of order divisible by $p$ are nilpotent.

Considering the coprime action of groups, as another generalization of Schmidt's theorem, Beltr\'{a}n and Shao [2, Theorem A] had the following result.

\begin{theorem} {\rm[2, Theorem A]}\ \  Let $G$ and $A$ be groups of coprime orders and assume that
$A$ acts on $G$ by automorphisms. If every maximal $A$-invariant subgroup of $G$ is nilpotent
but $G$ is not, then $G$ is solvable and $|G|=p^aq^b$ for two distinct primes $p$ and $q$, and $G$
has a normal $A$-invariant Sylow subgroup.
\end{theorem}

Moreover, Beltr\'{a}n and Shao [3, Theorem B] obtained the following result.

\begin{theorem} {\rm[3, Theorem B]}\ \ Suppose that a group $A$ acts coprimely on a group $G$ and
let $p$ be a prime divisor of the order of $G$. If every maximal $A$-invariant subgroup of $G$
whose order is divisible by $p$ is nilpotent, then $G$ is soluble.
\end{theorem}

In this paper, our main goal is to give a further generalization of above research by
considering the nilpotency of maximal invariant subgroup containing the normalizer of
some invariant Sylow subgroup, and we have the following result whose proof is given in
Section 3.

\begin{theorem}\ \ \label{thxx} Let $A$ and $G$ be groups such that $A$ acts coprimely on $G$ by automorphisms.
Then every maximal $A$-invariant subgroup of $G$ containing the normalizer of some
$A$-invariant Sylow subgroup is nilpotent if and only if one of following statements holds:

$(1)$ $G$ is nilpotent;

$(2)$ $G=(P_1\times\cdots\times P_{s-1})\times(P_s \rtimes(Q_1\times\cdots\times Q_t))$, where $P_1,\cdots,P_s$ are normal Sylow
subgroups of $G$ for $s\geq 1$, $Q_1,\cdots,Q_t$ are non-normal $A$-invariant Sylow subgroups of $G$
for $t\geq 1$, and there exists an $A$-invariant subgroup $E$ of $P_s$ such that $E$ is normal in $G$
and $E(Q_1\times\cdots\times Q_t)$ is a nilpotent maximal $A$-invariant subgroup of $P_s(Q_1\times\cdots\times Q_t)$.
\end{theorem}

By Theorem 1.3, the following corollary emerges.

\begin{cor}\ \ \label{c1} Every maximal subgroup of a group G containing the normalizer of some
Sylow subgroup is nilpotent if and only if one of following statements holds:

$(1)$ $G$ is nilpotent;

$(2)$ $G=(P_1\times\cdots\times P_{s-1})\times(P_s \rtimes(Q_1\times\cdots\times Q_t))$, where $P_1,\cdots,P_s$ are normal Sylow
subgroups of $G$ for $s\geq 1$, $Q_1,\cdots,Q_t$ are non-normal Sylow subgroups of $G$
for $t\geq 1$, and there exists a subgroup $E$ of $P_s$ such that $E$ is normal in $G$
and $E(Q_1\times\cdots\times Q_t)$ is a nilpotent maximal subgroup of $P_s(Q_1\times\cdots\times Q_t)$.
\end{cor}

\begin{remark}\ \ {\rm Let $A=1$ and $G=Z_5\times(Z_3\rtimes Z_2)$. This example shows that even if
every maximal $A$-invariant subgroup of $G$ containing the normalizer of some $A$-invariant
Sylow subgroup is nilpotent, we cannot get that every maximal $A$-invariant subgroup of
$G$ is nilpotent.}
\end{remark}

On non-nilpotent maximal subgroups, Shi [8, Theorem 5] proved that a group $G$ in
which every non-nilpotent maximal subgroup is normal has a Sylow tower. Furthermore,
Shi, Li and Shen [9, Theorem 1.7] showed that a group $G$ in which every non-nilpotent
maximal subgroup is normal is either $q$-nilpotent or $q$-closed for each prime divisor $q$ of
$|G|$. Considering the coprime action of groups, Li et al [5, Theorem 1.5] had the following
result.

\begin{theorem} {\rm[5, Theorem 1.5]}\ \ Suppose that $A$ acts on $G$ via automorphisms and that
$(|A|,|G|)=1$. If every non-nilpotent maximal $A$-invariant subgroup of $G$ is normal in $G$,
then $G$ is a Sylow tower group.
\end{theorem}

Moreover, Beltr\'{a}n and Shao [3, Theorem A] proved the following result.

\begin{theorem} {\rm[3, Theorem A]}\ \ Suppose that a group $A$ acts coprimely on a group $G$. If
every maximal $A$-invariant subgroup of $G$ is nilpotent or normal in $G$, then $G$ is $p$-nilpotent
for some prime $p$ dividing the order of $G$.
\end{theorem}

As a generalization of [5, Theorem 1.5] and [3, Theorem A], Shi and Liu obtained the
following result.

\begin{theorem} {\rm[10, Theorem 1.5]}\ \  Let $A$ and $G$ be groups such that $A$ acts coprimely on $G$
by automorphisms. If every non-nilpotent maximal $A$-invariant subgroup of $G$ is normal,
then $G$ is $p$-nilpotent or $p$-closed for each prime divisor $p$ of $|G|$.
\end{theorem}

Following from Theorem 1.3, we have the following result which provides a complete
extension of [5, Theorem 1.5], [3, Theorem A] and [10, Theorem 1.5].

\begin{theorem}\ \ \label{thxy} Let $G$ be a non-nilpotent group and $A$ be a group acting coprimely on $G$
by automorphisms. Then the following three statements are equivalent.

$(1)$ Every non-nilpotent maximal $A$-invariant subgroup of $G$ is normal;

$(2)$ Every maximal $A$-invariant subgroup of $G$ containing the normalizer of some $A$-invariant Sylow subgroup is nilpotent;

$(3)$ $G=(P_1\times\cdots\times P_{s-1})\times(P_s \rtimes(Q_1\times\cdots\times Q_t))$, where $P_1,\cdots,P_s$ are normal Sylow
subgroups of $G$ for $s\geq 1$, $Q_1,\cdots,Q_t$ are non-normal $A$-invariant Sylow subgroups of $G$
for $t\geq 1$, and there exists an $A$-invariant subgroup $E$ of $P_s$ such that $E$ is normal in $G$
and $E(Q_1\times\cdots\times Q_t)$ is a nilpotent maximal $A$-invariant subgroup of $P_s(Q_1\times\cdots\times Q_t)$.
\end{theorem}

The following corollary is a direct consequence of Theorem 1.9.

\begin{cor}\ \ \label{c2} Let $G$ be a non-nilpotent group. Then every non-nilpotent maximal
subgroup of $G$ is normal if and only if $G=(P_1\times\cdots\times P_{s-1})\times(P_s \rtimes(Q_1\times\cdots\times Q_t))$, where $P_1,\cdots,P_s$ are normal Sylow
subgroups of $G$ for $s\geq 1$, $Q_1,\cdots,Q_t$ are non-normal Sylow subgroups of $G$
for $t\geq 1$, and there exists a subgroup $E$ of $P_s$ such that $E$ is normal in $G$
and $E(Q_1\times\cdots\times Q_t)$ is a nilpotent maximal subgroup of $P_s(Q_1\times\cdots\times Q_t)$.
\end{cor}

Note that Corollary 1.10 provides a complete improvement of [9, Theorem 1.7].

In [10, Theorem 1.8] Shi and Liu obtained the following result.

\begin{theorem} {\rm[10, Theorem 1.8]}\ \ Let $A$ and $G$ be groups such that $A$ acts coprimely
on $G$ by automorphisms. Then every non-nilpotent maximal $A$-invariant subgroup of $G$
is a TI-subgroup if and only if every non-nilpotent maximal $A$-invariant subgroup of $G$ is
normal.
\end{theorem}

The following corollary is a direct consequence of Theorem 1.9 and [10, Theorem 1.8].

\begin{cor}\ \ \label{c3} Let $G$ be a non-nilpotent group and $A$ be a group acting coprimely on $G$
by automorphisms. Then every non-nilpotent maximal $A$-invariant subgroup of $G$ is a TI-subgroup if and only if
$G=(P_1\times\cdots\times P_{s-1})\times(P_s \rtimes(Q_1\times\cdots\times Q_t))$, where $P_1,\cdots,P_s$ are normal Sylow
subgroups of $G$ for $s\geq 1$, $Q_1,\cdots,Q_t$ are non-normal $A$-invariant Sylow subgroups of $G$
for $t\geq 1$, and there exists an $A$-invariant subgroup $E$ of $P_s$ such that $E$ is normal in $G$
and $E(Q_1\times\cdots\times Q_t)$ is a nilpotent maximal $A$-invariant subgroup of $P_s(Q_1\times\cdots\times Q_t)$.
\end{cor}

\section{Some Necessary Lemmas}\label{s2}

\begin{lemma} {\rm[1, Lemma 2.3]}\ \ Let $A$ be a group acting coprimely on a $p$-solvable group $G$.
If $H$ is a maximal $A$-invariant subgroup of $G$, then $|G:H|$ is a $p$-number or a $p'$-number.
\end{lemma}

\begin{lemma} {\rm[2, Theorem B]}\ \ Let $G$ and $A$ be groups of coprime orders and assume that
$A$ acts on $G$ by automorphisms. If $G$ has a nilpotent maximal $A$-invariant subgroup of odd
order, then $G$ is solvable.
\end{lemma}

\begin{lemma} {\rm[6, Lemma 2.3]}\ \ Let $A$ and $G$ be groups such that $A$ acts coprimely on $G$
by automorphisms. Suppose that $M$ is a maximal $A$-invariant subgroup of $G$, then $M$ is
either self-normalizing in $G$ or normal in $G$.
\end{lemma}

\begin{lemma} {\rm[12, Lemma 9]}\ \ Let $H$ be a nilpotent Hall-subgroup of a group $G$ that is not
a Sylow subgroup of $G$. If for each prime divisor $p$ of $|H|$, assume $P\in{\rm Syl}_p(H)$ we always
have $N_G(P)=H$, then there exists a normal subgroup $K$ of $G$ such that $G=KH$ and $K\cap H=1$.
\end{lemma}

\section{Proof of Theorem~1.3}\label{s3}

\demo {\bf Part I}. For the necessity part.

We first prove that $G$ is solvable.

Let $G$ be a counterexample of minimal order and $D\neq1$ be any $A$-invariant normal subgroup of $G$. It is easy to see that the hypothesis of the theorem also holds for the quotient group $G/D$. Then by the minimality of $G$, $G/D$ is solvable, which implies that $D$ is non-solvable.

Let $N$ be a minimal $A$-invariant normal subgroup of $G$ and $p$ be the largest prime divisor of $|N|$, then $p$ is odd since $N$ is non-solvable. Let $P$ be an $A$-invariant Sylow $p$-subgroup of $G$. Since $G$ has no non-trivial solvable $A$-invariant normal subgroups, $P$ is not normal in $G$. Then $N_G(P)<G$. It follows that there exists a maximal $A$-invariant subgroup $M$ of $G$ such that $N_G(P)<M$. By the hypothesis $M$ is nilpotent.

We claim that $M$ is a Hall-subgroup of $G$.

Otherwise, if $M$ is not a Hall-subgroup of $G$, then $M$ has an $A$-invariant Sylow $q$-subgroup $Q$ which is not a Sylow $q$-subgroup of $G$. Since $M$ is nilpotent, one has $N_G(Q)>M$. By the maximality of $M$, $N_G(Q) = G$, which implies that $Q$ is normal in $G$, a contradiction. Therefore, $M$ is a Hall-subgroup of $G$.

Since $p$ is odd and a group having a nilpotent maximal $A$-invariant subgroup of odd order is solvable by Lemma 2.2, $M$ is not a Sylow subgroup of $G$ by our assumption.

For any prime divisor $r$ of $|M|$, let $R$ be an $A$-invariant Sylow $r$-subgroup of $M$. Arguing as above, one must have $N_G(R) = M$.

So by Lemma 2.4, there exists a normal subgroup $H$ of $G$ such that $G=MH$ and $M\cap H=1$. Then $G=H\rtimes M$. Note that $H$ is a normal Hall-subgroup of $G$, then $H$ is a characteristic subgroup of $G$ which is also an $A$-invariant subgroup of $G$. Moreover, $H$ is a minimal $A$-invariant normal subgroup of $G$ by the maximality of $M$. Since $p$ is not dividing $|H|$, it follows that $N\neq H$. Then $N\cap H=1$. So $G\cong G/(N\cap H)$ which is isomorphic to a subgroup of $(G/N)\times(G/H)$. Thus $G$ is solvable since both $G/N$ and $G/H$ are solvable, a contradiction.

Hence the counterexample of minimal order does not exist and so $G$ is solvable.

Next we show that $G$ has normal Sylow subgroups.

Let $G$ be a counterexample of minimal order and $W$ be a minimal $A$-invariant normal subgroup of $G$. Since $G$ is solvable and $W$ is a characteristic simple group, $W$ is an elementary abelian group of order $w^{m}$, where $w$ is a prime divisor of $|G|$. By the minimality of $G$, $G/W$ has a normal Sylow subgroup $TW/W$, where $T$ is an $A$-invariant Sylow $t$-subgroup of $G$. By the Frattini argument, one has $G=N_G(T)TW=N_G(T)W$. Since $T$ is not normal in $G$ by our assumption, $N_G(T)<G$. There exists a maximal $A$-invariant subgroup $K$ of $G$ such that $N_G(T)<K$. Then $G=KW$. One has that $K$ is nilpotent by the hypothesis. Let $K_{w}$ be a Sylow $w$-subgroup of $K$, then $K_{w}W$ is a normal Sylow $w$-subgroup of $G$, a contradiction.

Hence the counterexample of minimal order does not exist and then $G$ has normal Sylow subgroups.

Final Conclusion.

Case (1): When $G$ is nilpotent, then the theorem is obviously true.

Case (2): When $G$ is non-nilpotent. Let $L=P_1\times P_2\times\cdots\times P_s$, where $P_1, P_2, \cdots,P_s$ are all normal Sylow subgroups of $G$. Since $G$ is non-nilpotent, $L<G$. For the quotient group $G/L$, arguing as above, $G/L$ has a normal Sylow subgroup $UL/L$, where $U$ is an $A$-invariant Sylow subgroup of $G$ and $U$ is not normal in $G$. Since $G$ is solvable, let $V$ be an $A$-invariant Hall-subgroup of $G$ such that $U\leq V$ and $G=L\rtimes V$. Note that $U=U(L\cap V)=UL\cap V$ is normal in $V$, then $V\leq N_G(U)$. Let $S$ be a maximal $A$-invariant subgroup of $G$ such that $N_G(U)\leq S$. By the hypothesis, $S$ is nilpotent, which
implies that $V$ is nilpotent.

Let $V=Q_1\times\cdots\times Q_t$, where $Q_1, \cdots, Q_t$ are non-normal $A$-invariant Sylow subgroups of $G$ for $t\geq1$. For every $1\leq i\leq t$, $V\leq N_G(Q_{i})$. Let $F$ be a maximal $A$-invariant subgroup of $G$ such that $N_G(Q_{i})\leq F$. By the hypothesis, $F$ is nilpotent. Since $G$ is solvable, $|G:D|$ is a prime-power by Lemma 2.1. Assume $P_1\times\cdots\times P_{s-1}\leq F$ and $P_s\nleq F$. Then $F=(P_1\times\cdots\times P_{s-1})\times(P_s\cap F)\times(Q_1\times\cdots\times Q_t)$. Let $E=P_s\cap F$, which is also an $A$-invariant subgroup of $G$. It is easy to see that $E<P_S$ and then $N_G(E)>F$. Therefore, $N_G(E)=G$ by the maximality of $F$. It follows that $E$ is normal in $G$. Since $F$ is a maximal $A$-invariant subgroup of $G$, one has that $E\times(Q_1\times\cdots\times Q_t)$ is a nilpotent
maximal $A$-invariant subgroup of $P_s\rtimes(Q_1\times\cdots\times Q_t)$.

{\bf Part II}. For the sufficiency part.

If group $G$ belongs to Case (1), the argument is trivial.

Next assume that group $G$ belongs to Case (2). For every non-normal $A$-invariant Sylow subgroup $Q_i$ of $G$, where $1\leq i\leq t$. By the hypothesis, it is easy to see that $N_G(Q_{i})=(P_1\times\cdots\times P_{s-1})\times E\times(Q_1\times\cdots\times Q_t)$ is a nilpotent maximal $A$-invariant subgroup of $G$.\hfill\qed

\section{Proof of Theorem~1.9}\label{s4}

\demo $(1)\Rightarrow(2)$ is trivial.

By Theorem 1.3 $(2)\Rightarrow(3)$ holds.

In the following we prove that $(3)\Rightarrow(1)$ holds.

Suppose that $M$ is any non-nilpotent maximal $A$-invariant subgroup of $G$. Let $H=P_1\times P_2\times\cdots\times P_s$ and $K=Q_1\times\cdots\times Q_t$. It is obvious that $G$ is solvable. Then $|G:M|$ is a prime-power by Lemma 2.1.

We divide our arguments into two cases.

Case $(i)$: Suppose $H<M$. Then $M=M\cap(H\rtimes K)=H\rtimes(M\cap K)$. By the maximality of $M$, $M\cap K$ is a maximal $A$-invariant subgroup of $K$. Since $K$ is nilpotent, $M\cap K$ is normal in $K$ by Lemma 2.3. It follows that $M$ is normal in $G$.

Case $(ii)$: Suppose $H\nleq M$. Then $G=HM$. Since $K$ is a Hall-subgroup of $G$ and $|G:M|$ is a prime-power, we can assume $K<M$. Therefore, $M=M\cap(H\rtimes K)=(H\cap M)\rtimes K$ and then $|H:H\cap M|=|G:M|$ is a prime-power.

If there exists a $P_i$ for $1\leq i\leq s-1$ such that $P_i\nleq M$, then $M=(P_1\times\cdots\times P_{i-1}\times(P_i\cap M)\times P_{i+1}\times\cdots P_{s-1})\times(P_s\rtimes K)$. It is obvious that $P_i\cap M$ is a maximal $A$-invariant subgroup of $P_i$ and then $P_i\cap M$ is normal in $P_i$. It follows that $M$ is normal in $G$.

Next assume $P_s\nleq M$. Then $M=(P_1\times\cdots\times P_{s-1})\times((P_s\cap M)\rtimes K)$.

If $E=1$, then $(P_1\times\cdots\times P_{s-1})\times K$ is a nilpotent maximal $A$-invariant subgroup of $G$ by the hypothesis, a contradiction.

Therefore, $E\neq1$. Moreover, $E\nleq P_s\cap M$. Note that $P_s\cap M$ is a maximal $A$-invariant subgroup of $P_s$. One has $P_s=E(P_s\cap M)$. Observe that $N_{P_{s}}(P_s\cap M)>P_s\cap M$, which implies that $N_{P_{s}\rtimes K}(P_s\cap M)>(P_s\cap M)\rtimes K$. Since $(P_s\cap M)\rtimes K$ is a maximal $A$-invariant subgroup of $P_s\rtimes K$, one has that $P_s\cap M$ is normal in $P_s\rtimes K$. Therefore, $P_s\cap M$ is normal in $G$. It follows that $M=(P_1\times\cdots\times P_{s-1})\times((P_s\cap M)\rtimes K)$ is normal in $E((P_1\times\cdots\times P_{s-1})\times((P_s\cap M)\rtimes K))=(P_1\times\cdots\times P_{s-1})\times((E(P_s\cap M))\rtimes K)=G$.\hfill\qed

\end{document}